\documentclass[draft]{amsart}
\usepackage{amsfonts, amsmath, amscd,amssymb,fullpage}
\usepackage{amssymb, amsbsy, amsthm,  amstext, amsopn, verbatim,
multicol}

\input{xy}

\xyoption{all}

\numberwithin{equation}{subsection}
\newtheorem*{thm}{Theorem}
\newtheorem*{prop}{Proposition}
\newtheorem*{lem}{Lemma}

 \newtheorem*{remark}{Remark}
\newtheorem*{remarks}{Remarks}

\newcommand{\im}{\operatorname{im}}

\newcommand{\und}[1]{\underline{#1}}

\newcommand{\ep}{\epsilon}

\newcommand{\g}{\mathfrak{g}}
\newcommand{\h}{\mathfrak{h}}

\newcommand{\Z}{\mathbb{Z}}

\newcommand{\C}{\mathbb{C}}

 \newcommand{\ttt}{\textsf}

\newcommand{\hr}{\mathfrak{h}^{\text{reg}}}

\newcommand{\EE}{\mathbf{E}}

\newcommand{\md}{\text{-mod}}
\DeclareMathOperator{\gl}{\mathfrak{gl}}

\DeclareMathOperator{\vect}{\ttt{Vect}}

\DeclareMathOperator{\id}{id}

\DeclareMathOperator{\gr}{gr}

\DeclareMathOperator{\pg}{\mathfrak{pg}}
\DeclareMathOperator{\rp}{Rep} \DeclareMathOperator{\mt}{Mat}
\DeclareMathOperator{\Tr}{Tr} \DeclareMathOperator{\diag}{diag}
 \DeclareMathOperator{\rank}{rank}
\DeclareMathOperator{\ed}{End} \DeclareMathOperator{\lie}{Lie}
\DeclareMathOperator{\pr}{\mathbb{P}}
  
\begin{document}

\title[Cherednik algebras and differential operators]{A remark on rational Cherednik algebras and differential operators on the cyclic quiver}
\begin{abstract}We show that the spherical subalgebra $U_{k,c}$ of the rational Cherednik algebra associated to $S_n\wr C_{\ell}$, the wreath product of the symmetric group and the cyclic group of order $\ell$, is isomorphic to a quotient of the ring of invariant differential operators on a space of representations of the cyclic quiver of size $\ell$. This confirms a version of \cite[Conjecture 11.22]{EG} in the case of cyclic groups. The proof is a straightforward application of work of Oblomkov, \cite{ob}, on the deformed Harish--Chandra homomorphism, and of Crawley--Boevey, \cite{CBgeom} and \cite{CBdecomp}, and Gan and Ginzburg, \cite{GG}, on preprojective algebras.
\end{abstract}

\maketitle
\section{Introduction}

\subsection{} The representation theory of symplectic reflection algebras has links with a number of subjects including algebraic combinatorics, resolutions of singularities, Lie theory and integrable systems. There is a family of symplectic reflection algebras associated to any symplectic vector space $V$ and finite subgroup $\Gamma \leq Sp(V)$, but a simple reduction allows one to study those subgroups $\Gamma$ which are generated by symplectic reflections (i.e. by elements whose set of fixed points is of codimension two in $V$). This essentially focuses attention on two cases:
\begin{enumerate}
\item $\Gamma = W$, a finite complex reflection group, acting on
$V=\h\oplus \h^*$ where $\h$ is a reflection representation of
$W$; \item $\Gamma = S_n \wr G$, where $G$ is a finite
subgroup of $SL_2(\C)$, acting naturally on $(\C^{2})^n$.
\end{enumerate}
 The representation theory in the first case is mysterious at the moment: several important results are known but there is no general theory yet. On the other hand a geometric point of view on the representation theory in the second case is beginning to emerge. A key fact is that in this case the singular space $V/\Gamma$ admits a crepant resolution of singularities: the representation theory of the symplectic reflection algebra is then expected to be closely related to the resolution. In the case $\Gamma = S_n$ (i.e. $G$ is trivial) there are two approaches to this: the first is via noncommutative algebraic geometry, \cite{GS}, the second via sheaves of differential operators, \cite{GG}. In this paper we extend the second approach to the groups $\Gamma = \Gamma_n = S_n\wr C_{\ell}$.

\subsection{} To state the result here we need to introduce a little notation. Let $Q$ be the cyclic quiver with $\ell$ vertices and cyclic orientation. Choose an extending vertex (in this case any vertex) $0$. Then let $Q_{\infty}$ be the quiver obtained by adding one vertex named $\infty$ to $Q$ that is joined to $0$ by a single arrow.

We will consider representation spaces of these quivers. Let
$\delta = (1,1,\ldots , 1)$ be the affine dimension vector of $Q$,
and set $\ep = e_{\infty} + n\delta$, a dimension vector for
$Q_{\infty}$. Let $\rp(Q, n\delta)$ and $\rp(Q_{\infty}, \ep)$ be the representation spaces of these quivers with the given dimension vectors. There is an action of $G = \prod_{r=0}^{\ell-1} GL_n(\C)$ on both these spaces. In fact, the action of the scalar matrices in $G$ is trivial on $\rp(Q,n\delta)$ (but not on $\rp(Q_{\infty}, \ep)$) so in this case the action descends to an action of $PG = G/\C^*$.

Let $\mathfrak{X} = \rp(Q,n\delta ) \times \mathbb{P}^{n-1}$. There is an action of $PG$ on $\mathfrak{X}$.

\subsection{}
Let $D(\rp(Q_{\infty}, \ep))$ denote the ring of differential operators on the affine space $\rp(Q_{\infty}, \ep)$, $D_{\mathfrak{X}}(nk)$ the sheaf of twisted differential operators on $\mathfrak{X}$ and $D(\mathfrak{X}, nk)$ its algebra of global sections. The group action of $G$ (respectively $PG$) on $\rp(Q_{\infty}, \ep)$ (respectively $\mathfrak{X}$) differentiates to an action of $\g = \lie(G)$ (respectively $\pg = \lie(PG)$) by differential operators. This gives mappings $$\hat{\tau}: \g \longrightarrow D(\rp(Q_{\infty}, \ep)), \qquad \tau : \pg \longrightarrow D_{\mathfrak{X}}(nk).$$

\subsection{}
\label{whatthm?}
Let $U_{k,c}$ be the spherical subalgebra of type $S_n\wr C_{\ell}$ (this is defined in Section \ref{spher}).
\begin{thm}
For all $(k,c)$ there are isomorphisms of algebras $$\left(\frac{D(\rp(Q_{\infty}, \ep))}{I_{k,c}}\right)^G \cong \left(\frac{D(\mathfrak{X} , nk)}{I_{c}}\right)^{PG} \cong U_{k,c},$$ where $I_{k,c}$ is the left ideal of $D(\rp(Q_{\infty}, \ep)$ generated by $(\hat{\tau} - \chi_{k,c})(\g)$ and $I_c$ is the left ideal of $D(\mathfrak{X}, nk)$ generated by $(\tau - \chi_c)(\pg)$ for suitable characters $\chi_{k,c}\in \g^*$ and $\chi_c\in \pg^*$ (which are defined in Section 4).
\end{thm}
Note that it is a standard fact that the left hand side is an
algebra. The proof of the theorem has two parts. One part
constructs a filtered homomorphism from the left hand side to the
right hand side using as its main input the work of Oblomkov,
\cite{ob}. The other part proves that the associated graded
homomorphism is an isomorphism and is a simple application of
results of Crawley--Boevey, \cite{CBgeom} and \cite{CBdecomp}, and
of Gan--Ginzburg, \cite{GG}.

\subsection{} We give an application of this result in Section 4.

\subsection{} While writing this down, we were informed that the general version of \cite[Conjecture 11.22]{EG} has been proved in \cite{EGGO}. That result is more general than the work presented here and requires a new approach and ideas to overcome problems that simply do not arise for the case $\Gamma = S_n \wr C_{\ell}$.

\section{Quivers}

\subsection{}
Once and for all fix integers $\ell$ and $n$. We assume that both
are greater than $1$. Set $\eta = \exp(2\pi i/\ell)$.

\subsection{}
Let $Q$ be the cyclic quiver with $\ell$ vertices and cyclic orientation. Choose an extending vertex (in this case any vertex) $0$. Then let $Q_{\infty}$ be the quiver obtained by adding one vertex named $\infty$ to $Q$ that is joined to $0$ by a single arrow. Let $\overline{Q}$ and $\overline{Q}_{\infty}$ denote the double quivers of $Q$ and $Q_{\infty}$ respectively.

We will consider representation spaces of these quivers. Let
$\delta = (1,1,\ldots , 1)$ be the affine dimension vector of $Q$,
and set $\ep = e_{\infty} + n\delta$, a dimension vector for
$Q_{\infty}$. Recall that $$\rp(Q, n\delta) =
\bigoplus_{r=0}^{\ell -1} \mt_n(\C) = \{ (X_0, X_1, \ldots , X_{\ell-1})\} = \{ (X) \}$$ and $$\rp(Q_{\infty}, \ep) =   \bigoplus_{r=0}^{\ell -1} \mt_n(\C)
\oplus \C^n = \{ (X_0, X_1,\ldots , X_{\ell-1}, i)\} = \{ (X, i)\}.$$ Let $G =
\prod_{r=0}^{\ell -1} GL_n(\C)$ be the base change group. If $g =
(g_0,\ldots , g_{\ell-1})$ then $g$ acts on $\rp(Q, n\delta)$ by
$$g\cdot (X_0, X_1, \ldots, X_{\ell-1}) = (g_0X_0g_1^{-1}, g_1X_1g_2^{-1} , \ldots , g_{\ell-1}X_{\ell-1}g_0^{-1})$$ and on $\rp(Q_{\infty}, \ep)$ by
$$g\cdot ( X_0, X_1 , \ldots , X_{\ell-1}, i) = ( g_0X_0g_1^{-1}, g_1X_1g_2^{-1} , \ldots , g_{\ell-1}X_{\ell-1}g_0^{-1}, g_0i).$$ The action of the scalar subgroup $\C^*$ is trivial in the first action (but not the second), so we can consider the first action as a $PG$--action where $PG = G/\C^*$. Let $\g$ and $\pg$ be the Lie algebras of $G$ and $PG$ respectively.

\subsection{}
Let $\hr \subset \C^n$ be the affine open subvariety consisting of
points $x = (x_1, \ldots , x_n)$ such that \begin{enumerate}
\item[(i)]
 if $i\neq j$ then $x_i \neq \eta^m x_j$  for all $m\in \Z$,  \item[(ii)] for each $1\leq i\leq n$ $x_i \neq 0$.\end{enumerate} This is
 the subset of $\C^n$ on which $\Gamma_n = S_n \wr C_{\ell}$ acts freely.

\subsection{} We can embed $\hr$ into $\rp(Q,n\delta)$ by first considering a point $x= (x_1, \ldots , x_n) \in \hr$ as a diagonal matrix $X = \diag(x_1,\ldots , x_n)$ and then sending this to $\und{X} = (X, X, \ldots ,X )$. We denote the image of $\hr$ in $\rp(Q,n\delta)$ by $\mathcal{S}$.

Let $T_{\Delta}$ be the subgroup of $G$ with elements $(T,T,\ldots
,T)$ where $T$ is a diagonal matrix in $GL_n(\C)$. Then
$T_{\Delta}$ is the stabiliser of $\mathcal{S}$. So consider the
mapping
$$\pi : G/T_{\Delta} \times \hr \longrightarrow \rp(Q,n\delta)$$
given by $\pi( gT_{\Delta}, x) = g\cdot \und{X}$. If we let $G$
act on $G/T_{\Delta}\times \hr$ by left multiplication then $\pi$
is a $G$--equivariant mapping.
\begin{lem}
\label{etale}
 $\pi$ is an \'{e}tale mapping with
covering group $\Gamma_n$. In fact its image
$\rp(Q,n\delta)^{reg}$ is open in $\rp(Q,n\delta)$ and we have an
isomorphism $$\omega: G/T_{\Delta} \times_{\Gamma_n} \hr
\longrightarrow \rp(Q,n\delta)^{reg} .$$
\end{lem}
\begin{proof}
 Let $\mathcal{S}=\{ \und{X} : x\in \hr\}$ and set $N_G(\mathcal{S})= \{ g\in G: g \cdot\mathcal{S} = \mathcal{S}\}$ and $Z_G(\mathcal{S}) = \{ g\in G: g\cdot \und X = \und X \text{ for all } \und X\in \mathcal{S}\}$.

 Suppose
$g\cdot \und{X} = \und{Y}$ for some $\und{X}, \und{Y}\in
\mathcal{S}$. This implies that for each $0\leq i\leq \ell-1$
$$g_i \diag(x)^{\ell} g_i^{-1} =
\diag(y)^{\ell}.$$ The hypotheses on $\hr$ imply that both
$\diag(x)^{\ell}$ and $\diag(y)^{\ell}$ are regular semisimple in
$\C^n$. Two such elements are conjugate if and only if $g_i \in
N_{GL_n(\C)}(T) = T\cdot S_n$ where $T$ is the diagonal subgroup
of $GL_n(\C)$. So there exists $\sigma \in S_n$ such that for all
$i$ we have $g_i = t_i \sigma$ for some $t_i \in T$, and for all
$1\leq r\leq n$ we have that $x_{\sigma (r)}^{\ell} =
y_{r}^{\ell}$. Hence $x_{\sigma(r)} = \eta^{m_r} y_{r}$ for some
$m_r\in \Z$. Now we find that $\und{Y} = g\cdot \und{X}$ implies
that $\diag(y_r) = t_it_{i+1}^{-1} \diag(\eta^{m_r}y_r) $. Since
$y_r \neq 0$ this shows that $t_{i+1} = \diag(\eta^{m_r}) t_{i}$
for each $i$. Hence we find that $gT_{\Delta} = (\sigma,
\diag(\eta^{m_r})\sigma, \ldots, \diag(\eta^{m_r})^{\ell-1}\sigma)
T_{\Delta}.$

In particular, if $\und{X} = \und{Y}$ we see from above that each
$m_r = 0$, so that $Z_{G}(\mathcal{S}) = T_{\Delta}$. Thus the
group $\Gamma_n$ is isomorphic to
$N_G(\mathcal{S})/Z_G(\mathcal{S})$ via the homomorphism that
sends $(\eta^{m_1},\ldots , \eta^{m_r})\sigma$ to $(\sigma,
\diag(\eta^{m_r})\sigma, \ldots,
\diag(\eta^{m_r})^{\ell-1}\sigma)T_{\Delta}$.

Now suppose that $\pi (gT_{\Delta}, x) = \pi (hT_{\Delta}, y)$.
Then $(h^{-1}g)\cdot \und X = \und Y$ and so we see that
$h^{-1}g\in N_G(\mathcal{S})$. This shows that $\pi$ is the
composition
$$ \begin{CD} G/T_{\Delta} \times \hr @>>> G/T_{\Delta} \times_{\Gamma_n} \hr @> \sim >> \rp(Q,n\delta)^{reg}\end{CD}.$$
The first mapping factors out the action of $\Gamma_n$, and since
$\Gamma_n$ acts freely on $\hr$ this is an \'{e}tale mapping.
Hence, to finish the lemma, it suffices to show that $\rp(Q,
n\delta)^{reg}$ is open in $\rp(Q,n\delta)$.

We claim first that $\rp(Q,n\delta)^{reg}$ is the set $O$ of
representations of $Q$ which decompose into $n$ simple modules of
dimension $\delta$ and whose endomorphism ring is $n$--dimensional.
To prove this observe that any element of $\rp(Q,n\delta)^{reg}$
is isomorphic to a representation of the form $\und{X}$ and so it decomposes
into the $n$ indecomposable modules $\und{X}_1, \ldots \und{X}_n$
of dimension $\delta$ where $\und{X}_i = (x_i, x_i, \ldots , x_i)$
(the condition $x_i\neq 0$ implies simplicity). Now the
representation $\und{X}_i$ is isomorphic to the representation
$(1,1,\ldots ,1,x_i^{\ell})$. By hypothesis $x_i^{\ell}\neq
x_j^{\ell}$ so we deduce that the representations $\und{x}_i$ are
pairwise non--isomorphic which ensures that the endomorphism ring
of $\und{X}$ is $n$--dimensional. This proves the inclusion $\rp (Q,n\delta)^{reg} \subseteq O$. On the
other hand, if $V$ belongs to $O$ then $V = V_1\oplus \ldots \oplus V_n$
where each $V_i$ is isomorphic to a representation $(1,1,\ldots
,1, \nu_i)$ for some non--zero scalars $\nu_i$. Moreover, since
$\dim \ed(V) = n$ the $\nu_i$ must be pairwise distinct. Now, let
$\eta_i$ be an $\ell$--th root of $\nu_i$. Then $V_i$ is
isomorphic to $(\eta_i, \ldots , \eta_i)$. Therefore $V$ is
isomorphic to the representation $\und{X}$ where $x = (\eta_1,
\ldots, \eta_n)$.

Now we must show that $O$ is open in $\rp(Q,n\delta)$. We use
first the fact that the canonical decomposition of the vector
$n\delta$ is $\delta + \delta + \cdots + \delta$, \cite[Theorem
3.6]{sc}. This means that the representations of $\rp(Q,n\delta)$
whose indecomposable components all have dimension $\delta$ form
an open set. Now, consider the morphism $f$ from $\rp(Q,\delta)$
to $\C$ which sends the representation $(\lambda_1 ,\ldots ,
\lambda_{\ell})$ to the product $\lambda_1\ldots \lambda_{\ell}$.
The open set $f^{-1}(\C^*)$ consists of the simple representations
of dimension vector $\delta$. Therefore the subset of
$\rp(Q,n\delta)$ consisting of representations which decompose as
the sum of $n$ simple representations of dimension vector $\delta$
is open. On the other hand, the function from $\rp(Q,n\delta)$ to
$\mathbb{N}$ which sends a representation $V$ to $\dim \ed(V)$ is
upper semi--continuous. Thus $\{ V : \dim \ed(v) \leq n\}$ is an
open set in $\rp(Q,n\delta)$. Intersecting these two sets shows
that $O$ is open, as required.
\end{proof}

\subsection{}
Now we're going to move from $Q$ to $Q_{\infty}$. So let's start
with the following $$\{ ([gT_{\Delta}, x], i) : g_0^{-1}i \text{
is a cyclic vector for }\diag(x) \} \subset
(G/T_{\Delta}\times_{\Gamma_n} \hr)\times \C^n.$$ By applying
$\omega^{-1}\times \id_{\C^n}$ this corresponds to an open subset
of $\rp(Q, n\delta) \times \C^n = \rp(Q_{\infty}, \ep).$ Call that
set $U_{\infty}$. This is a $G$--invariant open set since the
$G$--action on triples is given by $$h\cdot ([gT_{\Delta}, x], i)
= ([hgT_{\Delta}, x], h_0i)$$ so $g_0^{-1}i$ is cyclic for
$\diag(x)$ if and only if $(h_0g_0)^{-1}h_0 i$ is cyclic for
$\diag(x)$. Observe too that $U_{\infty}$ is an affine variety.
Indeed it is defined by the non--vanishing of the morphism $$s :
(G/T_{\Delta}\times_{\Gamma_n} \hr)\times \C^n \longrightarrow
\C$$ which sends $([gT_{\Delta}, x],i)$ to $(g_0^{-1}i) \wedge
\diag(x)\cdot(g_0^{-1}i) \wedge \cdots \wedge
\diag(x)^{n-1}\cdot(g_0^{-1}i)$.

\begin{lem}
\label{prinGbdle} The $G$--action on $U_{\infty}$ is free and
projection
 onto the second component $$\pi_2 : U_{\infty} \longrightarrow \hr/\Gamma_n$$
  is a principal $G$--bundle.
\end{lem}

\begin{proof}
Suppose that $h\cdot ([gT_{\Delta}, x], i) = ([gT_{\Delta}, x],
i)$.Then $[g^{-1}hgT_{\Delta}, x] =[T_{\Delta}, x]$, so by Lemma
\ref{etale} $g^{-1}hg \in T_{\Delta}$.

We have that $h_0i = i$. Setting $i' = g_0^{-1}i$ implies that
$g_0^{-1}h_0 g_0 i' = i'$. By hypothesis $i'$ is a cyclic vector
for $\diag(x)$. So in the standard basis $i'$ decomposes as $\sum
\lambda_je_j$ where each $\lambda_j$ is non--zero. Therefore the
only diagonal matrix that fixes $i'$ is the identity element. In
other words $g_0^{-1}h_0g_0 = I_n$. Since $g^{-1}hg\in T_{\Delta}$
this implies that $g^{-1}hg = \id$. Thus $h= \id$ and this proves
that the action is free.

It remains to prove that each fibre of $\pi_2$ is a $G$--orbit. So
take $([gT_{\Delta}, x], i) \in \pi_2^{-1}([x])$. This equals
$g\cdot ([T_{\Delta},x], g_0^{-1}i)$. Now $g_0^{-1}i$ is a cyclic
vector for $\diag(x)$ so it has the form $\sum \lambda_j e_j$ with
each $\lambda_j$ non--zero. Let $t= \diag(\lambda_1, \ldots,
\lambda_n)$ and consider $\und{t}= (t,\ldots, t)\in T_{\Delta}$.
We have $$([gT_{\Delta}, x], i) = g\und{t} \und{t}^{-1}
([T_{\Delta}, x], g_0^{-1}i) = g\und{t} ([T_{\Delta}, x] ,
\sum_{j=1}^n e_j).$$ This proves that each fibre of $\pi_2$ is
indeed a $G$--orbit.
\end{proof}

\subsection{}
Let $\rp(\overline{Q}_{\infty}, \ep)$ be the representation space
for the doubled quiver $\overline{Q}_{\infty}$: we can naturally
identify it with $T^*\rp(Q_{\infty}, \ep)$. The group $G$ acts on
the base and hence on the total space of the cotangent bundle. The
resulting moment map $$\mu : \rp(\overline{Q}_{\infty}, \ep)
\longrightarrow \g^*\cong \g$$ is given by $$\mu(X,Y,i,j) = [X,Y]
+ ij.$$

\begin{thm}[Gan--Ginzburg, Crawley--Boevey] \label{red} Let $\mu^{-1}(0)$ denote the scheme--theoretic fibre of $\mu$.
\begin{enumerate}
\item $\mu^{-1}(0)$ is reduced, equidimensional and a complete
intersection. \item The moment map $\mu$ is flat. \item
$\C[\mu^{-1}(0)]^{G} \cong \C[\h \oplus \h^*]^{\Gamma_n}.$
\end{enumerate}
\end{thm}

\begin{proof}

(i) This is \cite[Theorem 3.2.3]{GG}.

(ii) This follows from \cite[Theorem 1.1]{CBgeom} and the
dimension formula in \cite[Theorem 3.2.3(iii)]{GG}.

(iii) This is \cite[Theorem 1.1]{CBdecomp}
\end{proof}
\subsection{}
\label{xdef} Let $\mathfrak{X} = \{(X, i) \in \rp(Q,n\delta )
\times \mathbb{P}^{n-1}\}$. This space is the quotient of the
(quasi--affine) open subvariety $$U = \{ (X,i) : i\neq 0 \}
\subset \rp(Q_{\infty}, \ep)$$ by the scalar group $\C^*$. Thus
there is an action of $PG$ on $\mathfrak{X}$.

 Since
$$T^*\pr^{n-1} = \{ (i,j) : i\neq 0\text{, } ji = 0 \}/ \C^*$$ we have
$$T^*\mathfrak{X} = \{ (X, Y, i,j)\in \rp(\overline{Q}_{\infty},
\ep) : i\neq 0\text{, }ji=0\} /\C^*.$$ The $PG$ action on
$\mathfrak{X}$ gives rise to a moment map
$$\mu_{\mathfrak{X}} : T^*\mathfrak{X} \longrightarrow \pg^*\cong \pg.$$ Let
 $$\mu_{\mathfrak{X}}^{-1}(0) = \{ (X, Y, i,j)\in \rp(\overline{Q}_{\infty},
\ep) : i\neq 0\text{, }ji=0, \quad [X,Y]+ij =0\} / \C^*$$ denote
the scheme theoretic fibre of $0$.

\begin{prop}
\label{commprop} There is an isomorphism
$\C[\mu_{\mathfrak{X}}^{-1}(0)]^{PG} \cong \C[\h\oplus
\h^*]^{\Gamma_n}$. \end{prop}

\begin{proof}
Consider the $G$--equivariant open subvariety of
$\mu^{-1}(0)$ given by the non--vanishing of $i$. The variety
$\mu^{-1}(0)$ is determined by the conditions $[X,Y]+ ij =0,$ so
if we take the trace of this equation then we see that $0= Tr(ij)
= Tr(ji) = ji$. Thus we see that $\{ (X, Y, i,j)\in
\rp(\overline{Q}_{\infty}, \ep) : i\neq 0\text{, }ji=0\} \cap
\mu^{-1}(0)$ is an open subvariety of $\mu^{-1}(0)$ so in
particular reduced by Theorem \ref{red}(1). Hence factoring out by the
action of $\C^* \leq G$ shows that $\mu_{\mathfrak{X}}^{-1}(0)$ is
reduced and that there is a $PG$--equivariant morphism
$$\mu_{\mathfrak{X}}^{-1}(0) \longrightarrow \mu^{-1}(0)//\C^*.$$
This induces an algebra map
$$\alpha: \C[\mu^{-1}(0)]^{G}\longrightarrow
\C[\mu_{\mathfrak{X}}^{-1}(0)]^{PG}.$$

We now follow some of the proof of \cite[Lemma 6.3.2]{GG}. Write
$\ttt{O}_1$ for the conjugacy class of rank one nilpotent matrices
in $\gl(n)$, and let $\overline{\ttt{O}}_1$ denote the closure of
$\ttt{O}_1$ in $\gl(n)$. The moment map $\upsilon: T^*\pr^{n-1}
\longrightarrow \gl(n)^*\cong \gl(n)$ that sends $(i,j)$ to $ij$
gives a birational isomorphism $T^*\pr^{n-1} \longrightarrow
\overline{\ttt{O}}_1$. Let $J\subset \C[\gl(n)]=C[Z]$ be the ideal
generated by all $2\times 2$ minors of the matrix $Z$ and also by
the trace function. Then $J$ is a prime ideal whose zero scheme is
$\overline{\ttt{O}}_1$ and the pullback morphism $\upsilon^*:
\C[\gl(n)]/J \longrightarrow \C[T^*\pr^{n-1}]$ is a graded
isomorphism.

Now the moment map $\mu_{\mathfrak{X}}: T^*\mathfrak{X}
\longrightarrow \g^*$ factors as the composite
$$ \begin{CD} T^*\mathfrak{X} = T^*\rp(Q,n\delta) \times T^*\pr^{n-1}
@>>> T^*\rp(Q,n\delta) \times \overline{\ttt{O}}_1 @>{\theta} >>
\pg^*\end{CD}$$ where the first mapping is $\id \times \upsilon$
and the second mapping $\theta$ sends $(X,Y,Z)$ to $[X,Y] + Z_0$
where $Z_0$ indicates that we place the matrix $Z$ on the copy of
$\gl(n)$ associated to vertex $0$. We have a graded algebra
isomorphism
$$ \C[T^*\rp(Q,n\delta)]\otimes \C[\gl(n)]/J \longrightarrow
\C[T^*\mathcal{X}].$$ Now write $\C[X,Y,Z] = \C[T^*\rp(Q,n\delta)\
\times \gl(n)]$, and let $\C[X,Y,Z]([X,Y]+Z_0)$ denote the ideal
in $\C[X,Y,Z]$ generated by all matrix entries of the $\ell$
matrices $[X,Y]+Z_0$. Let $\mathbf{I}$ denote the ideal
$\C[X,Y,Z]([X,Y]+Z_0) + \C[X,Y] \otimes J \subset \C[X,Y,Z]$. From
the above we have
$$\C[\mu_{\mathfrak{X}}^{-1}(0)] \cong \C[T^*\rp(Q,n\delta) \times
\overline{\ttt{O}}_1]/\C[T^*\rp(Q,n\delta) \times
\overline{\ttt{O}}_1]\theta^*(\gl(n)) = \C[X,Y,Z]/\mathbf{I}.$$

Define an algebra homomorphism $r: \C[X,Y,Z] \longrightarrow
\C[X,Y]$ by sending $P\in \C[X,Y,Z]$ to the function $(X,Y) \mapsto P(X,Y,
-[X,Y]_0)$. Obviously $r$ induces an isomorphism
$\C[X,Y,Z]/\C[X,Y,Z]([X,Y]+Z_0) \cong \C[X,Y]/I_1$ where $I_1$ is
the ideal of $\C[\rp(\overline{Q},n\delta)]=\C[X,Y]$ generated by
the elements
$$ \sum_{h(a)=i} X_aX_{a^*} - \sum_{t(a) = i} X_{a^*}X_a$$ for
all $i$ not equal to zero. Observe that the linear function $P:
(X,Y,Z)\mapsto Tr Z = Tr ([X,Y]+Z_0)$ belongs to the ideal
$\C[X,Y,Z]([X,Y]+Z_0)$. We deduce that the mapping $r$ sends
$\C[X,Y]\otimes J$ to the ideal generated by $$\rank
(\sum_{h(a)=0} X_aX_{a^*} - \sum_{t(a) = 0} X_{a^*}X_a) \leq 1.$$
Thus we obtain algebra isomorphisms
$$\C[\mu_{\mathfrak{X}}^{-1}(0)] \cong \C[X,Y,Z]/\mathbf{I} \cong
\C[T^*\rp(Q,n\delta)]/I_2$$ where $I_2$ is ideal generated by the
elements
$$ \sum_{h(a)=i} X_aX_{a^*} - \sum_{t(a) = i} X_{a^*}X_a$$ for all
$1\leq i\leq \ell-1$, and $$\rank (\sum_{h(a)=0} X_aX_{a^*} -
\sum_{t(a) = 0} X_{a^*}X_a) \leq 1.$$

By \cite[Theorem 1]{lebpro} the $G$--invariant (respectively
$PG$--invariant) elements of $\C[\rp(\overline{Q}_{\infty}, \ep)]$
(respectively $\C[\rp(\overline{Q}, n\delta)]$) are generated by
traces along oriented cycles. Since all oriented cycles in
$\overline{Q}$ are oriented cycles in $\overline{Q}_{\infty}$ we have a
surjective composition of algebra homomorphisms \begin{equation}
\label{Xmom} \C[\h\oplus \h^*]^{\Gamma_n} \cong \C[\mu^{-1}(0)]^G
\longrightarrow \C[\mu_{\mathfrak{X}}^{-1}(0)]^{PG}
\longrightarrow
\left(\frac{\C[\rp(\overline{Q},n\delta)]}{I_2}\right)^{PG},\end{equation} where the first isomorphism is Theorem \ref{red}(3).
The left hand side is a domain of dimension $2\dim \h$, so to see
that the mapping is an isomorphism it suffices to prove that the
right hand side also has dimension $2\dim \h$.

 Let $I_3$ be the ideal of
$\C[\rp(\overline{Q},n\delta)]$ generated by the elements
$$ \sum_{h(a)=i} X_aX_{a^*} - \sum_{t(a) = i} X_{a^*}X_a$$ for
all $i$. This is the ideal of the zero fibre of the moment map for
the $PG$--action on $\rp(\overline{Q},n\delta)$. This ideal
contains $I_2$ since the rank condition on the matrices is implied
by the commutator condition. So there is a surjective mapping
$$\frac{\C[\rp(\overline{Q},n\delta)]^{PG}}{I_2^{PG}}
\longrightarrow
\frac{\C[\rp(\overline{Q},n\delta)]^{PG}}{I_3^{PG}}.$$ We do not
know yet whether the right hand side is reduced or not, but by
\cite[Theorem 1.1]{CBdecomp} the reduced quotient of the right
hand side is the ring of functions of the variety $(\h\oplus
\h^*)/\Gamma_n$. As this variety has dimension $2\dim \h$ we deduce that the composition in (\ref{Xmom}) is
an isomorphism, and hence that
$$\C[\mu_{\mathfrak{X}}^{-1}(0)]^{PG} \cong \C[\h\oplus
\h^*]^{\Gamma_n}.$$
\end{proof}

\begin{remark}
In passing let us note that the commutativity of the following
diagram
$$
\xymatrix{ \C[T^*\rp(Q,n\delta)] \ar@{->}[r]^-{\iota}
\ar@/_4pc/[drr]_{pr} & \C[T^*\rp(Q,n\delta)]\otimes \C[T^*\pr^{n-1}]
\ar@{->}[r] \ar@{->}[d]^{\upsilon^*} &
\C[\mu_{\mathfrak{X}}^{-1}(0)]
\ar@{->}[d]^{\wr} \\
& \C[T^*\rp(Q,n\delta)]\otimes \C[\overline{\ttt{O}}_1] \ar@{->}[r]^r  & \C[T^*\rp(q,n\delta)]/I_2 \\   &}
$$
where $\iota(f) = f\otimes 1$, shows that $\im \iota$ maps
surjectively onto $\C[\mu_{\mathfrak{X}}^{-1}(0)]$.
\end{remark}

\section{Differential operators}
\subsection{Symplectic reflection algebras} Let $C_{\ell}$ be the cyclic subgroup of $SL_2(\C)$ generated by $\sigma = \diag (\eta, \eta^{-1})$. The vector space $V = (\C^2)^n$ admits an action of $S_n\wr C_{\ell} = S_n \rtimes (C_{\ell})^n$: $(C_{\ell})^n$ acts by extending the natural action of $C_{\ell}$ on $\C^2$, whilst $S_n$ acts by permuting the $n$ copies of $\C^2$. For an element $\gamma\in C_{\ell}$ and an integer $1\leq i\leq n$ we write $\gamma_i$ to indicate the element $(1,\ldots, \gamma, \ldots ,1) \in C_{\ell}^n$ which is non--trivial in the $i$--th factor.

\subsection{}
The elements $S_n\wr C_{\ell}$ whose fixed points are a subspace of codimension two in $V$ are called symplectic reflections. In this case their conjugacy classes are of two types:
\begin{enumerate}
\item[($S$)] the elements $s_{ij}\gamma_i\gamma_j^{-1}$ where $1\leq i,j \leq n$, $s_{ij}\in S_n$ is the transposition that swaps $i$ and $j$, and $\gamma\in C_{\ell}$.

\item[($C_{\ell}$)] the elements $\gamma_i$ for $1\leq i \leq n$ and $\gamma \in C_{\ell}\setminus\{1\}$.
\end{enumerate}
There is a unique conjugacy class of type ($S$) and $\ell-1$ of type ($C_{\ell}$) (depending on the non--trivial element we choose from $C_{\ell}$). We will consider a conjugation invariant function from the set of symplectic reflections to $\C$. We can  identify it with a pair $(k,c)$ where $k \in \C$ and $c$ is an $\ell-1$--tuple of complex numbers: the function sends elements from ($S$) to $k$ and the elements $(\sigma^m)_i$ to $c_m$.

\subsection{} There is a symplectic form on $V$ which is induced from $n$ copies of the  standard symplectic form $\omega$ on $\C^2$. If we pick a basis $\{ x, y\}$ for $\C^2$ such that $\omega (x,y) = 1$ then we can extend this naturally to a basis $\{ x_i , y_i : 1\leq i\leq n\}$ of $V$ such that the $x$'s and the $y$'s form Lagrangian subspaces and $\omega (x_i, y_j) = \delta_{ij}$. We let $TV$ denote the tensor algebra on $V$: with our choice of basis this is just the free algebra on generators $x_i, y_i$ for $1\leq i\leq n$. The symplectic reflection algebra $H_{k,c}$ associated to $S_n\wr C_{\ell}$ is the quotient of $TV \ast (S_n\wr C_{\ell})$ by the following relations:
\begin{eqnarray*}
&x_i x_j = x_j x_i, \qquad y_i y_j = y_jy_i \qquad &\text{for all }1\leq i,j\leq n\\
&y_i x_i - x_iy_i = 1 + \frac{k}{2} \sum_{j\neq i} \sum_{\gamma \in C_{\ell}} s_{ij} \gamma_i \gamma_j^{-1} + \sum_{\gamma\in C_{\ell}\setminus\{1\}} c_{\gamma} \gamma_i \qquad &\text{for } 1\leq i\leq n\\
&y_i x_j - x_jy_i = -\frac{k}{2} \sum_{m=0}^{\ell-1} \eta^{m}s_{ij}(\sigma^m)_i(\sigma^m)_j^{-1} \qquad &\text{for } i\neq j.
\end{eqnarray*}
(NB: my $k$ is $-k$ for Oblomkov.)

\subsection{The spherical algebra} \label{spher} The symmetrising idempotent of the group algebra $C(S_n\wr C_{\ell})$ is
$$ e= \frac{1}{|S_n\wr C_{\ell}|} \sum_{w\in S_n\wr C_{\ell}} w.$$
The subalgebra $eH_{k,c}e$ is denoted by $U_{k,c}$ and called the {\it spherical algebra}. It will be our main object of study.

\subsection{Rings of differential operators}
Recall the definition of $\mathfrak{X}$ from \ref{xdef}. Let
$D_{\mathfrak{X}}(nk)$ denote the sheaf of twisted differential
operators on $\mathfrak{X}$ and let $D(\mathfrak{X}, nk)$ be its
algebra of global sections. This is simply the tensor product
$D(\rp(Q,n\delta)) \otimes D_{\mathbb{P}^{n-1}}(nk)$. (The twisted
differential operators on $\mathbb{P}^{n-1}$ can be defined as
follows. Let $A_n = \C[x_1,\ldots ,x_n, \partial_1,\ldots
,\partial_n]$ be the $n$--th Weyl algebra. This is a graded
algebra with $\deg(x_i) = 1$ and $\deg(\partial_i) = -1$. The
degree zero component is the subring generated by the operators
$x_i\partial_j$ which, under the commutator, generate the Lie
algebra $\mathfrak{gl}(n)$. Call this subring $R$. Let $\EE =
\sum_{i=1}^n x_i\partial_i \in R$ be the Euler operator. Then
$D(\mathbb{P}^{n-1}, nk)$ is the quotient of $R$ by the two--sided
ideal generated by $\EE - nk$.)

The group action of $PG$ on $\mathfrak{X}$ differentiates to an
action of $\pg$ on $\mathfrak{X}$ by differential operators. This
gives a mapping \begin{equation} \label{taudef} \tau : \pg
\longrightarrow D_{\mathfrak{X}}(nk).\end{equation} (One way to
understand this is to start back with $U\subset \rp(Q_{\infty},
\ep)$ and look at the $G$ action on $U$. Differentiating the
$G$--action gives an action of $\g$ by differential operators on
$U$, $\hat{\tau} : \g \longrightarrow D_U$. Since $\C^*$ acts
trivially on $\rp(Q, n\delta)$ and by scaling on $i\in
\rp(Q_{\infty}, \ep)$ we find that $\hat{\tau} ( \id) = 1\otimes
\EE$ where $\id = (I_n,I_n, \ldots , I_n)\in \C\subset \g$. Thus
we get an action of $\pg$ on $(D_U/D_U(1\otimes \EE- nk))^{\C^*} =
D_{\mathfrak{X}}(nk)$.)

\subsection{}
Recall the Lie algebra $\g = \lie (G)$ and its quotient $\pg = \lie(PG)$ which is simply $\g/\C\cdot \id$ where $\id = (I_n, \ldots, I_n) \in \g$. Let $\chi_c :\g \longrightarrow \C$ send  an element $(X) = (X_0, \ldots ,
X_{\ell-1})\in \g$ to $$\chi_c (X) = \sum_{r=0}^{\ell-1} C_r
\Tr(X_r)$$ where $C_r = \ell^{-1}(1 - \sum_{m=1}^{\ell-1}
\eta^{mr}c_m)$ for $1\leq r\leq \ell-1$ and $C_0 = \ell^{-1}(
1-\ell - \sum_{m=1}^{\ell-1} c_m)$. Observe that $$\chi_c(\id) = \Tr(I_n)
\sum_{r=0}^{\ell-1} C_i = n
\sum_{r=0}^{\ell-1}\sum_{m=0}^{\ell-1} - \eta^{rm} c_m = 0.$$ In
particular $\chi_c$ is actually a character of $\pg$.

Let $\chi_k : \g\longrightarrow \C$ send an element $(X) = (X_0,\ldots, X_{\ell-1})$ to $\chi_k(c) = k\Tr (X_0)$.

We will be regularly using the character $\chi_{k,c}\in \g^*$ defined by $\chi_{k,c} = \chi_c + \chi_k$.

\subsection{} Let us recall Oblomkov's deformed Harish--Chandra homomorphism, \cite{ob}. By Lemma \ref{etale} $\mathcal{S} = \omega (\hr/\Gamma_n)$
is a subset of $\rp(Q,n\delta)^{reg}$ which is a slice for the
$PG$--action on $\rp(Q,n\delta)$. Let
$$W'_k=(y_1\ldots y_n)^{-k}\C_{(0)}[y_1^{\pm 1}, \ldots , y_n^{\pm 1}],$$
a space of multivalued functions on $(\C^*)^n$. The Lie algebra
$\g$ acts on $W'_k$ by projection onto its $0$--th summand
$\gl(n)$, and then by the natural action of $\gl(n)$ on
polynomials (so $E_{ij}$ acts as $y_i\partial/\partial y_j$). With
this action the identity matrix in $\gl(n)$ becomes the Euler
operator
 $\mathbf{E}$ which acts by multiplication by $-nk$. Thus we can make $W'_k$ a
 $\pg$-module by twisting $W'_k$ by the character $\chi_{k}$ since then $\id$ acts trivially.
 If we call this module $W_k$ then $W_k = W'_k \otimes \chi_{k}$. Now define $Fun'$
 to be the space of functions on $\rp(Q,n\delta)$ of the form $$f = \tilde{f} \prod_{i=0}^{\ell
-1} \det(X_i)^{r_i}$$ where $\tilde{f}$ is a rational function on
$\rp(Q, n\delta)^{reg}$ regular on $\mathcal{S}$, $r_i =
\sum_{j=0}^i C_j + \sigma$ and $\sigma = \ell^{-1}
\sum_{s=0}^{\ell-1} sC_s$. Then $(Fun'\otimes W_k)^{\pg}$ is a
space of $(\pg,\chi_c)$--semiinvariant functions defined on a
neighbourhood of $\mathcal{S}$ which take values in $W_k$. This
space is a free $\C[\hr]^{\Gamma_n}$--module of rank $1$, the
isomorphism being given by restriction to $\mathcal{S}$. (Note
that the determinant of an element of the form $(X, \ldots ,X)$ is
$\det(X)^{\sum r_i} = 1$ as $\sum r_i = 0$.) Any $\pg$--invariant
differential operator, $D$, acts on such a function, $f$. Oblomkov
defines his homomorphism to be the restriction of $D(f)$ to
$\mathcal{S}$.

\subsection{} We can view the above procedure in terms of $\rp(Q_\infty, \ep)$.
Thanks to Lemma \ref{prinGbdle} we use $\mathcal{S}_{\infty} = \mathcal{S} \times (1,\ldots, 1)
\in U_{\infty}$ as a slice for the $G$--action. The space
$\mathcal{S}\times (\C^*)^n$ is a closed subset of $U_{\infty}$
since the condition that $i$ be cyclic for $\diag(x_1,\ldots ,
x_n)$ is equivalent to $i\in (\C^*)^n$. Thus functions on a
neighbourhood of $\mathcal{S}_{\infty}$
 in $U_{\infty}$
 can be identified with functions from a neighbourhood of
 $\mathcal{S}$ taking values in functions on $(\C^*)^n$. In
 particular, we can consider elements on $(Fun'\otimes W_k)^{\pg}$ first as $(\g,\chi_{k,c})$--semiinvariant functions
 from a neighbourhood of $\mathcal{S}$ taking values in $W_k'$ and hence
  as $(\g, \chi_{k,c})$--semiinvariant functions on
 an open set in a neighbourhood of $\mathcal{S}_{\infty}$.
 We can apply any element of $D\in D(U_{\infty})^{\g}$ to these $(\g, \chi_{k,c})$--semiinvariant functions
 and then
  restrict to $S_{\infty}$ to get a homomorphism
  $$\mathfrak{F}_{k,c} : D(U_{\infty})^{\g} \longrightarrow
  D(\hr/\Gamma_n).$$

\subsection{}
\label{obinfty} Since $\rp(Q_{\infty},\ep) = \rp(Q,n\delta) \times
\C^n$ there is a mapping $$\mathfrak{G}: D(\rp(Q,n\delta))^{\pg}
\longrightarrow D(U_{\infty})^{\g}$$ which sends $D \in
D(\rp(Q,n\delta))^{\pg}$ to $(D\otimes 1)$. Oblomkov's
homomorphism is $\mathfrak{F}_{k,c}\circ \mathfrak{G}$.

\subsection{} Differentiating the $G$--action on $U_{\infty}$
gives a Lie algebra homomorphism $\hat{\tau}: \g \longrightarrow
\vect(U_{\infty})$ which we extend to an algebra map
$$\hat{\tau}: U(\g) \longrightarrow D(U_{\infty}).$$ By Lemma \ref{prinGbdle} $U_{\infty}$ is a
principle $G$--bundle over $\hr/\Gamma_n$, so (a generalisation
of) \cite[Corollary 4.5]{sch} shows that the kernel of $\mathfrak{F}_{k,c}$
is $( D(U_{\infty})(\hat{\tau} - \chi_{k,c})(\g))^{\g}$. Moreover, since the
finite group $\Gamma_n$ acts freely on $\hr$ we can identify
$D(\hr/\Gamma_n)$ with $D(\hr)^{\Gamma_n}$.

\subsection{} Recall that $$D_{\mathfrak{X}}(nk) \cong \left(\frac{D_U}{D_U(\hat{\tau}
- \chi_{k})(\C\cdot \id)}\right)^{\C^*}.$$ Hence we have
\begin{equation}
\label{reduct}
\left(\frac{D_U}{D_U(\hat{\tau}-\chi_{k,c})(\g)}\right)^{G} \cong
\left(\frac{D_\mathfrak{X}(nk)}{D_\mathfrak{X}(nk)(\tau -
\chi_{c})(\pg)}\right)^{PG}, \end{equation} where $U = \{(X,i) :
i\neq 0\} \subset \rp(Q_{\infty},n\delta)$ as in \ref{xdef}. Now we
consider the restriction mapping $D_U \longrightarrow
D(U_{\infty})$. Composing the global sections of the above
isomorphism with this restriction and the homomorphism
$\mathfrak{F}_{k,c}$ gives
$$\mathfrak{R}'_{k,c} : \left(\frac{D(\mathfrak{X},nk)}{D(\mathfrak{X},nk)(\tau -
\chi_{c})(\pg)}\right)^{PG} \longrightarrow D(\hr)^{\Gamma_n}.$$

\subsection{} Let $$\delta_{k,c}(x) =
\delta^{-k-1}\delta^{\sigma}_{\Gamma}$$ where $\delta =
\prod_{1\leq i < j\leq n} (x_i^{\ell}-x_j^{\ell})$ and
$\delta_{\Gamma} = \prod_{i=1}^n x_i$. Define a twisted version of
$\mathfrak{R}'_{k,c}$ above $$\mathfrak{R}_{k,c}(D) =
\delta_{k,c}^{-1} \circ \mathfrak{R}'_{k,c}(D) \circ
\delta_{k,c}$$ for any differential operator $D$.

\subsection{} \label{mainthm} Our main result is the following.

\begin{thm} For all values of $k$ and $c$, the homomorphism $\mathfrak{R}_{k,c}$ has image $\im \theta_{k,c}$.
In particular we have an isomorphism $$\begin{CD} \theta_{k,c}^{-1}\circ
\mathfrak{R}_{k,c} :
\left(\frac{D(\mathfrak{X},nk)}{D(\mathfrak{X},nk)(\tau -
\chi_{c})(\pg)}\right)^{\pg}@> \sim >> U_{k,c}. \end{CD}$$ \end{thm}

\begin{proof}
Let us abuse notation by writing $U_{k,c}$ for the image of
$U_{k,c}$ in $D(\hr)^{\Gamma_n}$ under $\theta_{k,c}$.

Since $\mathfrak{X} = \rp(Q,n\delta) \times \pr^{n-1}$ there is a
mapping $$D(\rp(Q,n\delta))^{PG} \longrightarrow
D(\mathfrak{X},nk)^{PG}\longrightarrow D(\hr)^{\Gamma_n}$$ which
sends $D \in D(\rp(Q,n\delta))^{PG}$ to
$\mathfrak{R}_{k,c}(D\otimes 1)$. Recall $\tau$ from
\eqref{taudef}. Since $\gr \tau = \mu_{\mathcal{X}}^*$ we have an
inclusion $\gr(D(\mathfrak{X},nk))\mu_{\mathcal{X}}^*(\pg)
\subseteq \gr(D(\mathfrak{X},nk)(\tau-\chi_{c})(\pg))$. This gives
a graded surjection
$$p: \left(
\frac{\gr D(\mathfrak{X},nk)}{\gr(
D(\mathfrak{X},nk))\mu_{\mathcal{X}}^*(\pg)}\right)^{PG}\longrightarrow
\gr\left( \frac{D(\mathfrak{X},nk)}{D(\mathfrak{X},nk)(\tau -
\chi_{c})(\pg)}\right)^{PG}.$$ By Remark \ref{commprop} the
composition $$\gr D(\rp(Q,n\delta))^{PG}\longrightarrow \gr
D(\mathfrak{X},nk)^{PG} \longrightarrow \left(\frac{\gr
D(\mathfrak{X},nk)}{\gr(
D(\mathfrak{X},nk))\mu_{\mathcal{X}}^*(\pg)}\right)^{PG}\longrightarrow
\gr\left( \frac{D(\mathfrak{X},nk)}{D(\mathfrak{X},nk)(\tau -
\chi_{c})(\pg)}\right)^{PG}$$ is surjective. Thus the homomorphism
$$D(\rp(Q,n\delta))^{PG}\longrightarrow
\left(\frac{D(\mathfrak{X},nk)}{D(\mathfrak{X},nk)(\tau -
\chi_{c})(\pg)}\right)^{PG}$$ is also surjective. In particular,
by \ref{obinfty} this implies that the image of
$\mathfrak{R}_{k,c}$ equals the image of Oblomkov's
Harish--Chandra homomorphism, which, by \cite[Theorem 2.5]{ob}, is
$U_{k,c}$.

Thus we have a filtered surjective homomorphism
$$\mathfrak{R}_{k,c} :\left(\frac{D(\mathfrak{X},nk)}{D(\mathfrak{X},nk)(\tau -
\chi_{c})(\pg)}\right)^{PG} \longrightarrow U_{k,c}.$$ Thus the
dimension of the left hand side is at least $2\dim \h = \dim
U_{k,c}$. By Proposition \ref{commprop}
$$ \left(
\frac{\gr D(\mathfrak{X},nk)}{\gr(
D(\mathfrak{X},nk))\mu_{\mathcal{X}}^*(\pg)}\right)^{PG} \cong
\C[\mu_{\mathfrak{X}}^{-1}(0)]^{PG} \cong \C[\h\oplus
\h^*]^{\Gamma_n}.$$ Hence $p$ is a surjection from a domain of
dimension $2\dim \h$ onto an algebra of dimension at least $2\dim
\h$ and is hence an isomorphism. Thus
$\left(D(\mathfrak{X},nk)/D(\mathfrak{X},nk)(\tau -
\chi_{c})(\pg)\right)^{\pg}$ is a domain of dimension $2\dim \h$.
This implies that $\mathfrak{R}_{k,c}$ is an isomorphism.
\end{proof}

\section{Application: Shift functors}

\subsection{The Holland-Schwarz Lemma}
\label{holland} We want to understand the space
$$\frac{D(\rp(Q_{\infty},\ep))}{D(\rp(Q_\infty,\ep))(\hat{\tau} -
\chi_{k,c})(\g)}.$$ As we observed in the proof of Theorem
\ref{mainthm} there is a natural surjective homomorphism
\begin{equation} \label{hol} \frac{\gr D(\rp(Q_{\infty},\ep))}{\gr
D(\rp(Q_\infty,\ep))\mu^*(\g)}\longrightarrow
\gr\left(\frac{D(\rp(Q_{\infty},\ep))}{D(\rp(Q_\infty,\ep))(\hat{\tau}
- \chi_{k,c})(\g)}\right).\end{equation} It turns out that this is
an isomorphism.
\begin{lem}[Schwarz, Holland] The homomorphism \eqref{hol} is
an isomorphism of $\C[T^*\rp(Q_{\infty},\ep)]$--modules.
\end{lem}
\begin{proof} This is \cite[Lemma 2.2]{holl} since, by Theorem \ref{red}(2), the
moment map $\mu$ is flat.
\end{proof}
\subsection{} \label{bigmodule} This lets us prove the second part of the isomorphism in the statement of Theorem \ref{whatthm?}.
\begin{lem}
There is an algebra isomorphism
$$\left(\frac{D(\rp(Q_{\infty},\ep))}{D(\rp(Q_\infty,\ep))(\hat{\tau} -
\chi_{k,c})(\g)}\right)^{G}\longrightarrow
\left(\frac{D(\mathfrak{X},nk)}{D(\mathfrak{X},nk)(\tau -
\chi_{c})(\pg)}\right)^{PG}$$ \end{lem}
\begin{proof}
We have a natural $\pg$--equivariant mapping
$$D(\rp(Q_{\infty},\ep))^{\C^*} \longrightarrow D_U^{\C^*}
\longrightarrow D_{\mathfrak{X}}(nk)$$ which induces a
homomorphism
$$D(\rp(Q_{\infty},\ep))^{G} \longrightarrow \left(\frac{D(\mathfrak{X},nk)}{D(\mathfrak{X},nk)(\tau -
\chi_{c})(PG)}\right)^{\pg}.$$ This is surjective since, as we
observed in the proof of Theorem \ref{mainthm}, the image of
$D(\rp(Q,n\delta))^{PG} \subset D(\rp(Q_{\infty},\ep)^{G}$ spans
the right hand side. By \eqref{reduct} the kernel of this
homomorphism includes the ideal $(D(\rp(Q,\infty),\ep)(\hat{\tau}
- \chi_{k,c})(\g))^{G}$. Hence we have a surjective homomorphism
\begin{equation} \label{unproj} \left(\frac{D(\rp(Q_{\infty},\ep)}{D(\rp(Q,\infty),\ep)(\hat{\tau}
- \chi_{k,c})(\g)}\right)^{G}\longrightarrow
\left(\frac{D(\mathfrak{X},nk)}{D(\mathfrak{X},nk)(\tau -
\chi_{c})(\pg)}\right)^{PG}.\end{equation} By Lemma \ref{holland}
and Proposition \ref{commprop} there is an isomorphism
$$\left(\gr \frac{D(\rp(Q_{\infty},\ep)}{D(\rp(Q,\infty),\ep)(\hat{\tau} -
\chi_{k,c})(\g)}\right)^{G} \cong \left(\frac{\gr
D(\rp(Q_{\infty},\ep))}{\gr
D(\rp(Q_\infty,\ep))\mu^*(\g)}\right)^{G} = \C[\mu^{-1}(0)]^G =
\C[\h\oplus\h^*]^{\Gamma_n}.$$ This shows that the algebra on the
left is a domain of dimension of $2\dim \h$ and so \eqref{unproj}
is also injective, as required.
\end{proof}
\subsection{Shifting}
The previous two lemmas provide us with an interesting series of
bimodules. Given a character $\Lambda$ of $G$ we define
$$B_{k,c}^{\Lambda} = \left(\frac{D(\rp(Q_{\infty},\ep))}{D(\rp(Q_\infty,\ep))(\hat{\tau} -
\chi_{k,c})(\g)}\right)^{\Lambda}$$ to be the set of
$(G,\Lambda)$--semiinvariants. Thanks to Lemma \ref{bigmodule} and Theorem \ref{mainthm} this is a
right $U_{k,c}$--module. Now observe that if $x\in \g$ and $D\in
D(\rp(Q_{\infty},\ep))^{\Lambda}$ then
$$[\tau(x), D] = \lambda (x) D$$ where $\lambda = d\Lambda$. It follows that $B_{k,c}^{\Lambda}$ is also a left
$(D(\rp(Q_{\infty},\ep))/D(\rp(Q_\infty,\ep))(\hat{\tau} -
\chi_{k,c}-\lambda)(\g))^{G}$-module. So tensoring sets up a {\it shift functor}
$$ S_{k,c}^{\Lambda} : \left(\frac{D(\rp(Q_{\infty},\ep))}{D(\rp(Q_\infty,\ep))(\hat{\tau} -
\chi_{k,c})(\g)}\right)^{G}\md \longrightarrow
\left(\frac{D(\rp(Q_{\infty},\ep))}{D(\rp(Q_\infty,\ep))(\hat{\tau} -
\chi_{k,c}-\lambda )(\g)}\right)^{G}\md.$$

\subsection{} The character group of $G$ is
isomorphic to $\mathbb{Z}^{\ell}$ via $$(i_0,\ldots, i_{\ell-1})
\mapsto ((g_0,\ldots ,g_{\ell-1}) \mapsto \prod_{r=0}^{\ell-1}
\det(g_r)^{i_r}).$$ Corresponding to the standard basis element
$\ep_i$ is the character $\chi_i$ of $\g$ which sends $X\in \g$ to
$\Tr(X_i)$.
\label{param}
\begin{lem}
The bimodule corresponding to $\chi_i$ is a $(U_{k,c},
U_{k',c'})$--bimodule where $k'= k+1$ and $c' = c +
(1-\eta^{-i}, 1-\eta^{-2i},\ldots ,1- \eta^{-(\ell -1)i})$.
\end{lem}

\begin{proof}
Recall that $(k,c)$ corresponds to the character of $\g$ we called
$\chi_{k, c}$ which is defined as $$\chi_{k,c} (X) =
(C_0+k)\Tr(X_0) + \sum_{j=1}^{\ell-1} C_j\Tr(X_j),$$ where $C_r =
\ell^{-1}(1 - \sum_{m=1}^{\ell-1} \eta^{mr}c_m)$ for $1\leq r\leq
\ell-1$ and $C_0 = \ell^{-1}( 1-\ell - \sum_{m=1}^{\ell-1} c_m)$.
We need to calculate $(k',c')$ so that $\chi_{k,c} + \chi_i =
\chi_{k',c'}$. So we have $$(\chi_{nk,c}+\chi_i) (X) = (C_0+k)\Tr(X_0) + \Tr(X_i)+
\sum_{ j=1}^{\ell-1} C_j\Tr(X_j) = (C_0'+k')\Tr(X_0) + \sum_{j=1}^{\ell-1}
C'_j\Tr(X_j).$$ Calculation shows that $k' = k+1$
and that if $i=0$ then $C_j' = C_j$ and otherwise $$ C'_j = C_j + \begin{cases} -1 \quad &\text{if } j=0 \\ 1 & \text{if } j=i \\ 0 \quad &\text{otherwise.} \end{cases}$$ These unpack to give $c_m' = c_m + 1 -\eta^{-mi}$.
\end{proof}
\subsection{Question}
Thus for each $0\leq i \leq \ell-1$ we have a {\it shift functor}
$$S_i : U_{k,c}\md \longrightarrow U_{k+1, c'}\md$$ where $c'$ is as
above. When is this an equivalence of categories?

\begin{remarks} \begin{enumerate} \item We have been able to prove this is an equivalence when $(k,c)$ can be reached from $(0,0)$ by shifting.
\item Shift functors are also constructed in \cite{BC} and
\cite{RV}. Hopefully they agree with the functors here. \end{enumerate}
\end{remarks}


\begin{thebibliography}{GGOR}

\bibitem[BC]{BC} Y.~Berest and O.~Chalykh, Quasi--invariants of
complex reflection groups, {\it in preparation}.

\bibitem[BEG]{BEG} Y.~Berest, P.~Etingof and V.~Ginzburg, Cherednik algebras and differential operators on quasi--invariants, {\it Duke Math. J.} {\bf 118}, 279--337.

\bibitem[CB1]{CBgeom} W.~Crawley--Boevey, Geometry of the moment
map for representations of quivers, {\it Compositio Math.}, {\bf
126}, (2001), 257--293.

\bibitem[CB2]{CBdecomp} W.~Crawley--Boevey, Decomposition of
Marsden--Weinstein reductions for representations of quivers, {\it
Compositio Math.} {\bf 130} (2002),  225--239.

\bibitem[EG]{EG}  P.~Etingof and V.~Ginzburg, Symplectic reflection algebras, Calogero-Moser
   space, and deformed Harish-Chandra homomorphism,
    {\it Invent.\ Math.}, {\bf 147} (2002), 243-348.

\bibitem[EGGO]{EGGO} P.~Etingof, W.L.~Gan, V.~Ginzburg and A.~Oblomkov, private communication.

\bibitem[GG]{GG} W.L.~Gan and V.~Ginzburg, Almost commuting variety,
$\mathcal{D}$--modules, and Cherednik algebras,
 {\it RT:0409262}, March 2005.

\bibitem[GS]{GS} I.~Gordon and J.T.~Stafford, Rational Cherednik algebras and Hilbert schemes I and II, {\it to appear in Adv.Math. and Duke Math. Jour.}

\bibitem[H]{holl} M.~Holland, Quantization of the Marsden--Weinstein reduction for extended Dynkin quivers, {\it Ann. scient. \'{E}c. Norm. Sup.}, (1999), 813--834.

\bibitem[LP]{lebpro} L.~Le Bruyn and C.~Procesi, Semisimple
representations of quivers, {\it Trans.Amer.Math.Soc}, {\bf 317},
(1990), 585--598.

\bibitem[LS]{levsta} T.~Levasseur and J.T.~Stafford, The kernel of a homomorphism of Harish--Chandra, {\it Ann. scient. \'{E}c. Norm. Sup.}, {\bf 29}, (1996), 385--397.

\bibitem[O]{ob} A.~Oblomkov, Deformed Harish--Chandra homomorphism
for the cyclic quiver, {\it RT:0504395}, April 2005.

\bibitem[Scho]{sc} A.~Schofield, General representations of
quivers, {\it Proc. London Math. Soc} {\bf 65} (1992), 46--64.

\bibitem[Schw]{sch} G.W.~Schwarz, Lifting differential operators from orbit spaces, {\it Ann. scient. \'{E}c. Norm. Sup.}, {\bf 28}, (1995), 253--306.

\bibitem[V]{RV} R.~Vale, Diagonal coinvariants for $\Z_m\wr S_n$,{\it
RT:0505416}, May 2005.
\end{thebibliography}
\end{document}